\documentclass{SCAE}
\numberwithin{equation}{section}

\def\qed{\hfill$\Box$\par}
\def\cal{\mathcal}
\def\gg{\frak g}
\def\a{\alpha}
\def\b{\beta}
\def\d{\delta}
\def\D{\Delta}
\def\UU{{U}}

\def\e{\eta}

\def\vp{\varepsilon}

\def\v{\varphi}

\def\sc{\scriptstyle}
\def\ssc{\scriptscriptstyle}

\def\rar{\rightarrow}\def\rar{\to}

\def\ra{\rangle}
\def\AA{{\cal A}}

\def\la{\langle}
\def\ni{\noindent}
\def\ptl{\partial}

\def\N{\mathbb{N}{\ssc\,}}
\def\Z{\mathbb{Z}{\ssc\,}}

\def\F{\mathbb{F}{\ssc\,}}
\def\QED{\hfill$\Box$} \numberwithin{equation}{section}
\newtheorem{theo}{Theorem}[section]
\newtheorem{conv}[theo]{Convention}
\newtheorem{defi}[theo]{Definition}
\newtheorem{rema}[theo]{Remark}

\newtheorem{prop}[theo]{Proposition}

\def\adddot{$\!\!\!${\bf.}\ \ }
\def\adddot{}

\renewcommand{\baselinestretch}{1.3}
\begin{document}

\Year{2015} %
\Month{January}
\Vol{58} %
\No{1} %
\BeginPage{1} %
\EndPage{xx} %
\AuthorMark{Song G A, Su Y C }
\ReceivedDay{October 11, 2014}
\AcceptedDay{January 22, 2015}
\PublishedOnlineDay{; published online   2015}
\DOI{xx.xxxx/sxxxxx-xxx-xxxx-x} 

\title{Dual Lie Bialgebra Structures of  Poisson
Types}{}


\author[1]{ Guang'ai Song}{}
\author[2]{Yucai Su}{Corresponding author}

\address[{\rm1}]{College of Mathematics and Information Science, Shandong
Institute of Businessand Technology , Yantai {\rm 264005}, China;}
\address[{\rm2}]{Department of Mathematics, Tongji University,
Shanghai {\rm 200092},China}

\Emails{gasong@sdibt.edu.cn,
ycsu@tongji.edu.cn}\maketitle


 {\begin{center}
\parbox{14.5cm}{\begin{abstract}
 Let $\AA=\F[x,y]$ be the polynomial algebra on two variables $x,y$ over an algebraically closed field $\F$ of characteristic zero. Under the
Poisson bracket, $\AA$ is equipped with a natural Lie algebra structure.
    It is proven that  the maximal good subspace of $\AA^*$ induced from the multiplication of the associative commutative algebra $\AA$ coincides
 with the maximal good subspace of $\AA^*$ induced from the Poisson bracket of the Poisson Lie algebra $\AA$. Based on this,
structures of dual Lie bialgebras of the Poisson type are
investigated. As by-products, five classes of new infinite dimensional Lie algebras are obtained.
\end{abstract}}\end{center}}



\renewcommand{\baselinestretch}{1.2}
\begin{center} \renewcommand{\arraystretch}{1.5}
{\begin{tabular}{lp{0.8\textwidth}} \hline \scriptsize
{\bf Citation:}\!\!\!\!&\scriptsize Song G A, Su Y C. Science  China: Mathematics  title. Sci China Math, 2014, 57, doi: 10.1007/s11425-000-0000-0\vspace{1mm}
\\
\hline
\end{tabular}}\end{center}

\baselineskip 11pt\parindent=10.8pt  \wuhao
\section{Introduction}
Lie bialgebras, having
close relations with Yang-Baxter equations \cite{2}, are important ingredients in quantum groups, which have drawn more and more attentions in literature (e.g., \cite{1,1-2,D,G,M1,3,4,5,6,7,8,9,10,12,13,14}).
Michaelis  \cite{4} investigated  structures of Witt type Lie bialgebras.  Ng and Taft \cite{9} gave a classification of this
type Lie bialgebras, and obtained that all structures of Lie bialgebras  on the one sided
 Witt algebra, the Witt algebra and the Virasoro algebra are
 coboundary triangular (cf.~\cite{8}).
  For the cases of generalized Witt type Lie algebras and generalized Virasoro-like Lie algebras, the authors of \cite{10, 13} proved that
 all  structures of Lie bialgebras on them are
 coboundary triangular. Similar results hold
 for some other kinds of Lie algebras (cf., e.g., \cite{13, 14}).

From  the examples of infinite dimensional Lie bialgebras constructed in \cite{4}, many infinite
dimensional Lie bialgebra structures we know are coboundary
triangular. It may sound that coboundary triangular Lie bialgebras
are relatively simple. However, they are not trivial in the sense that many
natural problems associated with them remain open (see, also Remark \ref{Raa22}). For
example, even for  the (two-sided) Witt algebra and the Virasoro
algebra, a complete classification of coboundary triangular Lie bialgebra
structures on them is still an open problem. Nevertheless, not much on representations of infinite dimensional Lie bialgebras
is known. 
From the viewpoint of
Lie bialgebras, considering dual Lie bialgebra structures may help us
understand more on infinite dimensional Lie bialgebra structures. For instance, by considering  structures of  dual Lie
bialgebras of Witt and Virasoro types, the authors of \cite{16}
surprisingly obtained some new series of infinite
dimensional Lie algebras. 
In the present paper, we study structures of dual Lie bialgebras of Poisson
type.
 One may have noticed that the dual of a finite dimensional Lie
bialgebra is naturally a Lie bialgebra, and so one would not predict
anything new in this case. However, for the cases of infinite
dimensional Lie bialgebras, the situations become quite different, which can be seen in the following contents.

Let us recall the definition of Poisson algebras here: a   {\it Poisson algebra} is a triple $({\cal P},[\cdot,\cdot],\cdot)$ such that $({\cal P},[\cdot,\cdot])$ is a Lie algebra, $({\cal P},\cdot)$ is an associative algebra, and the following {\it Leibniz rule} holds:
\begin{equation}\label{Poi1}
[a,bc]=[a,b]c+b[a,c]\mbox{ \ for \ }a,b,c\in {\cal P}.
\end{equation}
In particular, for any commutative associative algebra $({\cal A},\cdot)$, and any commutative derivations $\ptl_1,\ptl_2$ of ${\cal A}$, we obtain a
{\it Poisson algebra}  $({\cal A},[\cdot,\cdot],\cdot)$ with Lie bracket $[\cdot,\cdot]$ defined as follows.
\begin{equation}\label{Poi2}
[a,b]=\ptl_1(a)\ptl_2(b)-\ptl_2(a)\ptl_1(b)\mbox{ \ for \ }a,b\in {\cal A}.
\end{equation}
If we take ${\cal A}=\F[x^{\pm1},y^{\pm1}]$ (where $\F$ is an algebraically closed  field of characteristic zero) and $\ptl_1\!=\!x\frac{\ptl}{\ptl x},\,\ptl_2\!=\!y\frac{\ptl}{\ptl y}$, then we obtain the {\it Virasoro-like algebra} $({\cal A},[\cdot,\cdot])$
with basis $\{x^iy^j\,|\,i,j\!\in\!\Z\}$ and Lie bracket
\begin{equation}\label{Poi3}
[x^iy^j,x^ky^\ell]=(i\ell-jk)x^{i+k}y^{j+\ell}\mbox{ \ for \ }i,j,k,\ell\in \Z.
\end{equation}
The Virasoro-like algebra \eqref{Poi3} can be generalized as follows:
For any {\it nondegenerate} additive subgroup $\Gamma$ of $\F^2$ (i.e.,  $\Gamma$ contains an $\F$-basis of $\F^2$), we have the group algebra ${\cal A}=\F[\Gamma]$ with basis $\{ L_{\a}\, |\, \a \in \Gamma \}$ and
 multiplication defined by $\mu(L_{\a}, L_{\b}) = L_{\a+\b}$ for $\a, \b \in \Gamma$.
Then we have the (generalized) Virasoro-like algebra $(\cal{A}, \v )$ with Lie bracket $\v$ defined by
  \begin{equation}\label{VL-1}
  \v(L_{\a}, L_{\b}) = (\a_1\b_2- \b_1\a_2)L_{\a+\b}
   \mbox{ \ for \ $\a=(\a_1, \a_2), \,\b=(\b_1, \b_2) \in \Gamma.$}
\end{equation}
Furthermore, if we take ${\cal A}=\F[x,y]$ and $\ptl_1=\frac{\ptl}{\ptl x},\,\ptl_2=\frac{\ptl}{\ptl y}$, then we obtain the {\it classical Poisson algebra} $(\F[x,y],[\cdot,\cdot],\cdot)$, whose
Lie bracket is given by
\begin{equation}\label{Poi4}
[f,g]=J(f,g)\mbox{ \ for \ }f,g\in\F[x,y],
\end{equation}
where $J(f,g):=\Big|{}^{\frac{\ptl f}{\ptl x}\ \frac{\ptl f}{\ptl y}}_{\frac{\ptl g}{\ptl x}\ \frac{\ptl g}{\ptl y}}\Big|$ is the {\it Jacobian determinant} of $f$ and $g$.

The reason  we have a special interest in the classical Poisson algebra also lies in the  fact that this algebra is closely related to the distinguished Jacobian conjecture (e.g., \cite{V2,Y}), which can be stated as ``any non-zero endomorphism of  $(\F[x,y],[\cdot,\cdot],\cdot)$ is an isomorphism''. One observes that a Jacobi pair $(f,g)$ (i.e., $f,g\in\F[x,y]$ satisfying $J(f,g)\in\F\backslash\{0\}$) corresponds to a solution $r=f\otimes fg-fg\otimes f$ of the classical Yang-Baxter Equation (cf.~\eqref{CYBE}), thus gives rise to a Lie bialgebra structure on $\F[x,y]$.

The paper is organized as follows. Some
 definitions and preliminary results are briefly recalled in Section 2. Then in Section 3, structures of dual coalgebras
 of $\F[x, y]$ are addressed. Finally in Section 4,  structures of
dual Lie bialgebras of Poisson type are investigated.
 The main results of the present paper are summarized in Theorems \ref{theo-3-1}, \ref{theorem-2}, \ref{th-2-2}, \ref{theorem-3-1}  and \ref{last-theo}.

\section{Definitions and preliminary
results}

Throughout the paper, 
all vector spaces are assumed to be over an algebraically closed field $\F$ of characteristic zero. As usual,
we use $\Z_+$ to denote the set of nonnegative
integers.
 We briefly recall some notions on Lie bialgebras, for details, we
refer readers to, e.g., \cite{2,10}.
\begin{defi}\adddot{\label{2.1}}\rm\begin{enumerate}\parskip-1pt\item A {\it Lie bialgebra}
is a triple $(L, [\cdot, \cdot], \d)$ such that
$(L, [\cdot, \cdot])$ is  a   Lie  algebra,
$(L, \d)$ is   a  Lie coalgebra, and $\d:L\to L\otimes L$ is a derivation, namely,
$\d[x, y] = x \cdot \d(y) - y\cdot\d(x)$ for  $x,y \in L$,
 where  $ x \cdot (y\otimes z) = [x, y]\otimes z + y
\otimes [x, z]$ for  $x, y, z \in L.$
\item A  Lie bialgebra
 $(L, [\cdot, \cdot],  \d)$ is {\it coboundary} if $\d$ is coboundary in the sense that there exists $r \in L \otimes L$ written as $r= \sum r^{[1]} \otimes r^{[2]} $, such that $\d(x) = x\cdot r$
for $x\in L$.
\item A coboundary Lie bialgebra
 $(L, [\cdot, \cdot],  \d)$ is {\it triangular} if $r$ satisfies  the following {\it classical Yang-Baxter Equation} (CYBE),
\begin{equation}\label{CYBE}C(r) = [r_{12}, r_{13}] + [r_{12}, r_{23}] + [r_{13}, r_{23}]
 =0,\end{equation}
 \noindent where
$r_{12} = \sum r^{[1]} \otimes r^{[2]} \otimes 1,$
$r_{13}= \sum r^{[1]} \otimes 1 \otimes r^{[2]},$
$r_{23} = \sum r^{[1]} \otimes 1 \otimes r^{[2]}$
are elements in $ \UU(L)\!
\otimes\! \UU(L) \!\otimes\! \UU(L)$, and $\UU(L)$ is the universal enveloping algebra of $L$.\end{enumerate}
 \end{defi}

Two Lie bialgebras $(\gg,[\cdot,\cdot],\d)$ and $(\gg',[\cdot,\cdot]',\d')$ are said to be {\it dually paired}
if their bialgebra structures are related via
\begin{equation}\label{Dual-paired}
\langle[f, h ]', \xi\rangle = \langle f \otimes h, \d \xi \rangle ,\ \ \  \langle \d' f, \xi \otimes
\eta \rangle = \langle f, [\xi, \eta]\rangle
\mbox{ \ for $f, h \in \gg',\ \xi, \eta \in \gg,$ }\end{equation}
 where $\langle\cdot,\cdot\rangle$ is a nondegenerate
bilinear form on $\gg'\times\gg$, which is naturally extended to a nondegenerate
bilinear form on \mbox{$(\gg'\otimes\gg')\times(\gg\otimes\gg)$.}
In particular, if $\gg'=\gg$ as a vector space, then $\gg$ is called a {\it
self-dual Lie bialgebra}.

The following  result whose proof is straightforward can be found in \cite{3}.
\begin{prop}\adddot
Let $(\gg, [\cdot, \cdot], \d)$ be a finite dimensional Lie
bialgebra, then so is the linear dual space $\gg ^*:={\rm Hom}_\F(\gg ,\F)$ by dualisation, namely $(\gg^*,[\cdot,\cdot]',\d')$ is the Lie bialgebra
defined by  \eqref{Dual-paired} with $\gg'=\gg^*$. In particular, $\gg$ and $\gg^*$ are dually paired.
\end{prop}

Thus a finite dimensional Lie biallgebra $(\gg,[\cdot,\cdot],\d)$ is always self-dual as there exists a vector space isomorphism $\gg\to\gg^*$ which pulls back the bialgebra structure on $\gg^*$ to $\gg$ to obtain another bialgebra structure on $\gg$ to make it to be self-dual.
However, in sharp contrast to the finite dimensional case, infinite
dimensional Lie bialgebras are  not self-dual in general.

For convenience, we denote by $\v$  the Lie bracket of Lie algebra $(\gg ,
[\cdot, \cdot])$, which can be regarded as a linear map $\v:\gg\otimes\gg\to\gg$. Let $\v^{\ast}:\gg^*\to(\gg\otimes\gg)^*$ be  the dual of $\v$.

\begin{defi}\adddot\rm \cite{M1}
Let $(\gg , \v)$ be a Lie algebra over $\F$. A  subspace $V$ of $\gg ^*$ is called a {\it good
subspace} if $\v^*(V) \subset V \otimes V.$
Denote $\Re = \{ V\, |\, V {\rm\ is \ a \ good \ subspace \ of \ }\gg ^*\}$. Then
$
\gg ^{\circ} = \mbox{$\sum
_{V\in \Re}$} V,$ 
 is also a good subspace
of $\gg ^*$, which is obviously the maximal good subspace
of $\gg ^*$.
\end{defi}

It is clear that
if $\gg $ is a finite dimensional Lie algebra, then
$\gg ^{\circ} = \gg ^*$.

\begin{prop}\adddot{\rm\cite{M1}}\label{PPPP} For any good subspace $V$ of $\gg ^*$, the pair $(V,
\v^{\ast})$ is a Lie coalgebra. In particular, $(\gg ^{\circ}, \v^*)$
is a Lie coalgebra.
\end{prop}

 For any Lie algebra $\gg $, the dual space $\gg ^*$ has a natural right $\gg $-module
 structure defined for $f\in \gg ^*$ and $ x\in \gg $ by $$(f\cdot x)(y) = f([x,
 y])\mbox{ \ for \ }y\in \gg .$$ We denote $f\cdot \gg  = {\rm span} \{ f\cdot x\, |\, x \in
 \gg \}$, the {\it space of translates} of $f$ by elements of $\gg $.

We summarize some results of \cite{1,1-2,D, G} as follows.
 \begin{prop}\adddot\label{Th111}
 Let $\gg $ be a Lie algebra.
 Then\begin{enumerate}\item $\gg ^{\circ}=\{f\in\gg^*\,|\,f\cdot\gg\mbox{ is finite dimensional}{\sc\,}\}.$
\item
$\gg ^{\circ}=(\v^*)^{-1}(\gg^*\otimes\gg^*)$, the preimage of $\gg^*\otimes\gg^*$ in $\gg^*$.
\end{enumerate}
 \end{prop}

The notion of good subspaces of an associative algebra can be defined analogously.
In the next two sections, we shall investigate $\gg^\circ$ for some associative or
Lie algebras $\gg$.

\section{The structure of $\F[x,y]^{\circ}$
}
Let $(\cal{A}, \mu, \eta)$ be an associative $\F$-algebra with unit,  where $\mu$ and $\eta$ are respectively the multiplication $\mu: \cal{A} \otimes
\cal{A} \rar \cal{A}$ and  the unit $\eta: \F \rar \cal{A}$,
satisfying $$\begin{array}{ll}
\mu \circ (id \otimes \mu)= \mu \circ (\mu \otimes id)):&
\cal{A} \otimes \cal{A} \otimes \cal{A} \rar \cal{A},\\[2pt]
(\eta\otimes id)(k\otimes a) = (id \otimes \eta)(a \otimes k): &\F \otimes
\cal{A}\cong \cal{A}\otimes \F \cong \cal{A},\end{array}$$
 for  $k \in \F,\, a \in \cal{A}$.
Then a coassociative coalgebra is a triple $(C, \D, \e)$,  which is obtained
by conversing arrows in the definition of an associative algebra.
Namely,  $\D: C \rar C\otimes C$ and $\e: \F \rar C$ are respectively
comultiplication and counit of $C$, satisfying
$$\begin{array}{ll}(\D \otimes id)\circ \D = (id \otimes \D)\circ \D: &C \rar
C\otimes C \otimes C,\\[2pt]
(\e \otimes id)\circ \D = (id \otimes
\e)\circ \D: &C \rar C\otimes C \rar \F\otimes C \cong C\otimes \F
\cong C.\end{array}$$
For any vector space $\cal{A}$, there exists a natural  injection $\rho:
\cal{A^{\ast}} \otimes \cal{A^{\ast}} \rar (\cal{A} \otimes
\cal{A})^{\ast}$ defined by $\rho(f, g )(a, b)= \langle f, a\rangle \langle g, b\rangle $ for
$f, g \in \cal{A^{ \ast}}$ and $ a, b, \in \cal{A}$. In case $\cal A$ is finite dimensional,
$\rho$ is an isomorphism.
If $(\cal{A}, \mu)$ is associative, the multiplication $\mu
$ induces the map $\mu^{\ast}: \cal{A^{\ast}}
\rar (\cal{A} \otimes \cal{A})^{\ast}$. If $\cal{A}$ is  finite
dimensional, then the isomorphism
$\rho$ insures that $(\cal{A^{\ast}}, \mu^{ \ast},
\e^{\ast})$ is a coalgebra, \vspace*{-4pt}where for simplicity, $\mu^{\ast}$ denotes the
composition of the maps: $ \cal{A^{\ast}}
\stackrel{\mu^{\ast}}{\rar}(\cal{A} \otimes \cal{A})^{\ast}
\stackrel{(\rho)^{{\rm -1}}}{\rar} \cal{A^{\ast}} \otimes
\cal{A^{\ast}}$.

%
%
%

Now let  $(\cal{A}, \mu)$ be a commutative associative algebra. 
Then $\cal{A^{\circ}} =
(\mu^{\ast})^{-1}(\cal{A^*} \otimes \cal{A^*})$ (cf.~\cite{12} and Proposition \ref{Th111}). For
$\ptl \in {\rm Der}(\cal{A})$ and $f \in \cal{A^{\circ}}$, using
$$\ptl \mu
=\mu( id \otimes \ptl + \ptl \otimes id),\ \ \ \mu^*\ptl^*(f) = (id
\otimes \ptl^* + \ptl^* \otimes id )\mu^*(f) \in \cal{A^*} \otimes
\cal{A^*},$$ we obtain $\ptl^*(\cal{A^{\circ}}) \subset
\cal{A^{\circ}}$. Thus, we observe that there are two natural approaches to produce  Lie
coalgebras from some subspaces of $\cal{A^*}$. One is induced
from the associative structure of $\cal{A}$ as
follows:  First we have the
cocommutative coassociative
coalgebra $(\cal{A^{\circ}}, \mu^{\circ})$ with $\mu^{\circ}:= \mu^*|_{\cal{A^{\circ}}}$. Then we obtion
the Lie coalgebra $\cal{A}^{\circ}_{\mu}:= (\cal{A^{\circ}}, \D)$ with cobracket, induced from
cocommutative coassociative coalgebra structure, defined by
\begin{equation}\label{co-b0} \D(f) = (\ptl_1^{\circ} \otimes \ptl_2^{\circ} -\ptl_2^{\circ} \otimes \ptl_1^{\circ} )
\mu^{\circ} (f) {\rm \ \ for\ }\ f \in \cal{A^{\circ}},
\end{equation}
where $\ptl_1, \ptl_2
\in {\rm Der}{\cal{A}}$ are two fixed derivations satisfying $\ptl_1\ptl_2= \ptl_2\ptl_1.$ Here and below, for  any $\ptl \in {\rm Der}(\cal{A}),$
we denote $\ptl^{\circ}=\ptl^*|_{\cal{A^{\circ}}}$.

Another approach is as follows: Let $\cal{A}_{\v} = (\cal{A}, [\cdot, \cdot])$ be the Lie algebra
defined in \eqref{Poi2} (where $\v = [\cdot, \cdot]$). The Lie coalgebra induced from $\cal{A}_{\v}$ is
$\cal{A}_{\v}^{\circ}= (\cal{A}_{\v}^{\circ}, \v^{\circ})$, where the subspace
$\AA_{\v}^{\circ}$ of $\AA^*$ is determined by Proposition \ref{Th111} with 
cobracket 
defined by
\begin{equation}\label{cob-1}  \v^{\circ}(f) = (\mu(\ptl_1\otimes \ptl_2 -
\ptl_2\otimes \ptl_1))^{\ast}(f)=(\ptl_1\otimes \ptl_2 -
\ptl_2\otimes \ptl_1)^{\ast}\mu^*(f){\rm \ \ for \ }\ f \in
\cal{A}_{\v}^{\circ}.
\end{equation}

\begin{prop}\adddot\label{p2} Let $(\cal{A}, \mu)$ be a commutative
associative algebra with unit, 
and $\ptl_1, \ptl_2 \in
{\rm Der}(\cal{A})$ are commutative. Then the Lie coalgebra $\cal{A}^{\circ}_{\mu}$ 
is a Lie subcoalgebra of $\cal{A}_{\v}^{\circ}$.
\end{prop}
\ni{\it Proof.~}~ For  $f \in \cal{A}^{\circ}_\mu$, we have
$\v^{\circ}(f)
=
(\ptl_1\otimes \ptl_2 - \ptl_2\otimes \ptl_1)^{\ast}\mu^*(f)
=
(\ptl_1\otimes \ptl_2 - \ptl_2\otimes
\ptl_1)^{\ast} \mu^{\circ}(f)
=(\ptl_1^{\circ}\otimes \ptl_2^{\circ} -
\ptl_2^{\circ}\ptl_1^{\circ})\mu^{\circ}(f)= \D(f),
$ 
where the last equality follows from \eqref{co-b0}. Thus, $\cal{A}_{\mu}^{\circ}$ is a Lie subcoalgebra of $\cal{A}^{\circ}_{\v}$.
\QED\vskip4pt

\begin{theo}\adddot\label{theo-3-1} Let $(\AA, \mu)$ be a commutative associative algebra, and $(\cal{A}, \v)$  the Poisson Lie  algebra
defined in \eqref{Poi2}. If there exists $h \in \AA$ such that the ideal $I$ of $(\AA,\v)$ generated by $h$ has finite codimension, then  $\cal{A}_{\mu}^{\circ} = \cal{A}_{\v}^{\circ}.$ In particularly, if $\AA= \F[x, y]$, then $\cal{A}_{\mu}^{\circ} = \cal{A}_{\v}^{\circ}.$
\end{theo}

\ni{\it Proof. } 
Denote $\cdot$ and $\star$ the actions of
 $(\cal{A}, \mu)$ and $(\cal{A}, \v)$ on $\cal{A}^*$ respectively, i.e.,
 $(f\cdot a)(b) = f(\mu(a, b))$, and $(f\star a)(b)= f(\v(a, b))$
for $a, b \in \cal{A}, f \in \cal{A}^*$. From the relation $\v(a, bc) = \v(a, c)b + \v(a, b)c$ for $ a, b, c \in \AA$, we have
\begin{equation*}
(f\star a)\cdot b - (f\cdot b)\star a = f\cdot \v(a, b), \ \forall\,a, b \in \AA.
\end{equation*}
If $f \in {\AA}_{\v}^{\circ}$, i.e. $f\star \AA$ is finite dimensional, then $f\cdot \v(b, \AA)$ is finite dimensional. Thus if the ideal $I$ has finite codimension, and $f\cdot I$ is finite dimensional, it follows that $f\cdot \AA$ is finite dimensional. From Proposition \ref{Th111}, we have $f\in \AA_{\mu}^{\circ}$.
 \QED

\begin{rema} \rm
The difference between  $\cal{A}_{\mu}^{\circ}$ and $\cal{A}_{\v}^{\circ}$ is that
${\cal{A}}_{\mu}^{\circ}$, as a Lie coalgebra, is induced from coassociative  coalgebra $(\AA^{\circ}, \mu^{\circ})$,
  and $\AA^{\circ}_{\mu}$ is determined by $(\mu^*)^{-1}(\AA^*\otimes \AA^*)$ as a vector subspace of $\AA^*$
  $($a good subspace of the dual of $(\AA, \mu)\,)$, but $\cal{A}_{\v}^{\circ}$ is the dual of the Lie algebra $(\AA, \v)$ (determined by Proposition \ref{Th111}).
\end{rema}
%

\section{Dual Lie bialgebras of Poisson type
}\setcounter{section}{4}\setcounter{theo}{0}\setcounter{equation}{0}
Poisson algebras (cf.~\eqref{Poi1}) have important algebra structures, which  have close relations with the Virasoro
algebra and vertex operator (super)algebras (e.g., \cite{AH}). They can be also regarded as special cases of Lie algebras of Block type.
Therefore, some attentions have been paid  on them and some related Lie algebras 
(e.g., \cite{CL,13,Su04,17,GGS,SXX13,SXX12,SXZ,SXY,11, X}).
In this section, 
we consider the dual structures of
Poisson type Lie bialgebras.

The following result can be found in \cite{12}.
\begin{prop} \adddot\label{S} Let $A, B$ be commutative associative algebras,
regarding $A^* \otimes B^* \subset (A \otimes B)^*$, then $A^{\circ}
\otimes B^{\circ} = (A \otimes B)^{\circ}.$
\end{prop}

Recall from \cite{10} that the dual space of $\F[x]$ can be identified with the space $\F[[\vp]]$. From \cite{7, 16}, and Proposition \ref{S}, we have

\begin{prop}\adddot\label{F[x]}\begin{enumerate}
\item Let $f= \sum_{i=0}^{\infty} f_i \vp^i \in
 \F[[\vp]]$ with $f_i\in\F$. Then
\begin{eqnarray*}
f \!\in\! \F[x]^{\circ}&\!\!\! \Longleftrightarrow\!\!\!& f_n \!=\! h_1f_{n-1}\! +\! h_2
f_{n-2}
\! +\! \cdots
\! +\! h_r f_{n-r}\mbox{ for some $r\in\N,\,h_i\!\in\!\F$ and all $n\!>\!r$}\\[-4pt]
&\!\!\!\Longleftrightarrow \!\!\!&f \in \Big\{\frac{g(\vp)}{h(\vp)}\,\Big|\,g(\vp),h(\vp)\in\F[\vp],\,h(0)\ne0\Big\}.
\end{eqnarray*}
 \item Denote $\cal{A}= \F[x, y],$ the polynomial algebra on tow variables $x,
 y$. Then
$$\cal{A}^{\circ} = \F[x, y]^{\circ} \cong \F[x]^{\circ}
 \otimes \F[y]^{\circ}.$$
 \end{enumerate}
\end{prop}
%

Let $(\gg , \v, \d)$  be a Lie bialgebra. The map 
$\v^*$ induces a map $\v^{\circ}:=\v^*|_{\gg^{\circ}}:\gg ^{\circ}
\rar \gg ^{\circ} \otimes \gg ^{\circ}$, making
$(\gg^{\circ},\v^\circ)$  to be a Lie coalgebra. By
\cite[Proposition 3]{8}, the map $\d^*:\gg ^* \otimes \gg ^*
\hookrightarrow (\gg \otimes \gg )^* \stackrel{\d^*}{\rar} \gg ^*$
induces a map
 $\d^{\circ}:=\d^*_{\gg ^{\circ} \otimes \gg ^{\circ}}:
\gg ^{\circ} \otimes \gg ^{\circ} \rar \gg ^{\circ}$, making
$(\gg^{\circ},\d^\circ)$ to be a Lie algebra. Thus we obtain a Lie bialgebra
  $(\gg ^{\circ}, \d^{\circ}, \v^{\circ} )$, the {dual Lie
 bialgebra} of $(\gg , \v,  \d)$.

\def\ell{l}Now take $\AA=\F[x, y]$. 
Let $\vp^i, \e^i$  be duals of $x^i, y^i \in \AA$ respectively,
namely, $\langle\vp^i\e^j,x^ky^\ell\rangle=\vp^i\e^j(x^ky^\ell)=\d_{i,k}\d_{j,\ell}$
for $i,j,k,\ell\in\Z_+$. 
Any 
element $u\in{\cal A}^*$ can be written as $u=\sum_{i,j}u_{i,j}\vp^i\e^j$ (possibly an infinite
sum). Let $g=\sum_{k,l}g_{k,l}x^ky^l \in{\cal A}$ (a finite sum).
Then
\begin{equation}\la u,g\ra=
u(g) = \mbox{$\sum\limits_{i,j, k,l}$}u_{i,j}g_{k,l}\langle \vp^i\e^j, x^ky^l
\rangle =\mbox{$\sum\limits_{i,j, k,l}$}u_{i,j}g_{k,l}\d_{i,k}\d_{j, l}
\mbox{ (a finite sum).}
\end{equation}
Let $\ptl_1=\frac{\partial}{\partial x},\, \ptl_2=\frac{\partial}{\partial y}$.
Then we have the Poisson Lie algebra $(\cal{A}, \v)$ defined by \eqref{Poi4}. 
From Theorem \ref{theo-3-1},
 it is easy to check that $\AA^{\circ}_{\v} = \AA^{\circ}_{\mu}$.

\begin{conv}\rm\label{Con1} \begin{itemize}\item[(1)]
 If an undefined notation appears in an expression, we treat it zero; for instance $\vp^i\e^j=0$  if $i<0\ \mbox{or}\ j<0.$
\item[(2)]
 When there is no confusion, we  use $[\cdot,\cdot]$ to denote the bracket in $\gg$ or  $\gg^\circ$, i.e., $[\cdot,\cdot]=\v$ or $\d^{\circ}$.
We also use $\D$ to denote the cobracket in $\gg$ or  $\gg^\circ$, i.e., $\D=\d$ or
 $\v^\circ$.
\end{itemize}\end{conv}

Let $m,n\in\Z_+$, and take $ a = x^my^n, \, b = xy \in \cal{A}$. 
Then $[a, b ] =
(m-n)a$. Thus by \cite{4}, the triple $(\cal{A}, [\cdot, \cdot],
\D_r)$ with $r = a\otimes b
- b \otimes a$ is a coboundary triangular Lie bialgebra whose  cobracket is
defined by
\begin{equation}\D_r (f) = f\cdot r= [f, a] \otimes b + a \otimes [f, b] - [f, b]\otimes a - b\otimes [f, a]\ {\rm
for \  }\ f \in \cal{A}.
\end{equation}
\begin{theo}\adddot \label{theorem-2}Let 
$(\AA, [\cdot, \cdot], \D_r)$ be
 the coboundary triangular Lie bialgebra defined above. 
The dual Lie bialgebra of $\AA$
is $(\AA^{\circ}, [\cdot, \cdot], \D)$, where $\AA^{\circ}$ is described by Proposition $\ref{F[x]}(2)$ with cobracket $\D$ defined by
\begin{equation} \label{dual-1}\D(\vp^m\e^n) =\mbox{$\sum\limits_{k+s=m+1,\ l+t=n+1}(kt-ls)$}\vp^k\e^l\otimes
\vp^s\e^t,
\end{equation}
and bracket $[\cdot,\cdot]$ uniquely determined by the skew-symmetry and the following
\begin{equation}\label{bracket-4-2}
[\vp^i\e^j,\vp^s\e^t]=\left\{
\begin{array}{ll}
(m(t\!+\!1)\!-\!n(s\!+\!1))\vp^{s+1-m}\e^{t+1-n}\!\!\!&\mbox{if \ }(i,j)=(1,1),\,(s,t) \neq (1,1),\\[4pt]
(s-t)\vp^s\e^t&\mbox{if \ }(i,j)=(m, n) \neq (s, t)\neq(1, 1), \\[4pt]
 0&\mbox{otherwise}.\end{array}\right.
\end{equation}
\end{theo}
\ni{\it Proof.~}~ 
Assume
$\mu^{\circ}
(\vp^m\e^n) = \mbox{$\sum_{k,l,.s,t \in
\Z_+}$}c_{k,l,s,t}\vp^k\e^l\otimes \vp^s\e^t\mbox{ for some }c_{k,l,s,t}
\in \F.$
\noindent Then
\begin{eqnarray*}\begin{array}{lll}
c_{i,j,p,q} =
\mu^{\circ}(\vp^m\e^n)(x^iy^j\otimes x^py^q)
=
\langle\vp^m\e^n,
\mu(x^iy^j\otimes x^py^q)\rangle
=
\langle\vp^m\e^n,
x^{i+p}y^{j+q}\rangle
=
\d_{m, i+p}\d_{n, j+q}.
\end{array}
\end{eqnarray*}
Thus, $
\mu^{\circ}(\vp^m\e^n) =\mbox{$\sum
_{k+s=m,\,
l+t=n}$}\vp^k\e^l\otimes \vp^s\e^t.
$ 
Assume $\ptl_1^{\circ}(\vp^i\e^j) =
\sum_{s,t}c_{s,t}\vp^s\e^t$. Then
$$c_{k,l} = 
\ptl_1^{\circ}(\vp^i\e^j)(x^ky^l)=
\vp^i\e^j(\ptl_1(x^ky^l))=k\d_{i,k-1}\d_{j,l}=(i+1)\d_{i+1, k}\d_{j, l},$$
\noindent i.e.,  $\ptl_1^{\circ}(\vp^i\e^j)
=(i+1)\vp^{i+1}\e^j.$ Similarly,  $\ptl_2^{\circ}(\vp^i\e^j)=
(j+1)\vp^i\e^{j+1}.$ From (\ref{co-b0}), we obtain
\begin{eqnarray*}\begin{array}{lll}
\D(\vp^m\e^n) &=& (\ptl_1^{\circ} \otimes \ptl_2^{\circ}-
\ptl_2^{\circ}\otimes\ptl_1^{\circ})\mu^{\circ}(\vp^m\e^j)
=
\sum\limits_{k+s=m+1,\, l+t=n+1}(kt-ls)\vp^k\e^l\otimes \vp^s\e^t.
\end{array}
\end{eqnarray*}
Therefore, (\ref{dual-1}) holds. Next, we verify \eqref{bracket-4-2}.
We have
\begin{eqnarray}
\!\!\!\!\!\!\!\!\!\!\!\!\!\!\!&\!\!\!\!\!\!\!\!\!\!\!\!\!\!\!\!\!\!\!\!\!\!\!\!\!\!\!\!\!\!&
\langle[\vp^i\e^j, \vp^s\e^t], x^ky^l \rangle
=\langle \vp^i\e^j\otimes \vp^s\e^t, (kn-lm)x^{k+m-1}y^{l+n-1}\otimes xy
+(k-l)x^my^n\otimes x^ky^l
\nonumber\\\!\!\!\!\!\!\!\!\!\!\!\!\!\!\!\!\!\!\!\!
&\!\!\!\!\!\!\!\!\!\!\!\!\!\!\!&
\phantom{\langle[\vp^i\e^j, \vp^s\e^t], x^ky^l \rangle=}
- (kn-lm)xy\otimes x^{k+m-1}y^{l+n-1}
-(k-l)x^ky^l\otimes x^my^n\rangle
\nonumber
\\ 
\label{f-3}\!\!\!\!\!\!\!\!\!\!\!\!\!\!\!\!\!\!\!\!
&\!\!\!\!\!\!\!\!\!\!\!\!\!\!\!&
\phantom{\langle[\vp^i\e^j, \vp^s\e^t], x^ky^l \rangle}=\langle(n(i+1)-m(j+1))\d_{s,1}\d_{t,1}\vp^{i+1-m}\e^{j+1-n}
+(s-t)\d_{i,m}\d_{j,n}\vp^s\e^t
\nonumber\\\!\!\!\!\!\!\!\!\!\!
&\!\!\!\!\!\!\!\!\!\!\!\!\!\!\!&
\phantom{\langle[\vp^i\e^j, \vp^s\e^t], x^ky^l \rangle=}
-(n(s\!+\!1)\!-\!m(t\!+\!1))\d_{i,1}\d_{j,1}\vp^{s+1-m}\e^{t+1-n}
\!-\!(i\!-\!j)\d_{s,m}\d_{t,n}\vp^i\e^j,x^ky^l\rangle.
\end{eqnarray}
If  $(i,j)=(1, 1)$ and $ (s, t)\neq (1, 1)$, then (\ref{f-3}) gives
(note that $(m, n)\neq (1, 1)$)
$$[\vp\e, \vp^s\e^t ] = \left\{\begin{array}{lll}(m(t+1)-n(s+1))\vp^{s+1-m}\e^{t+1-n}&\mbox{if}\  s+1-m\geq 0,\ t+1-m\geq 0,\\
0&\mbox{otherwise}.\end{array}\right.$$
Thus we obtain the first case  of (\ref{bracket-4-2}) (cf.~Convention \ref{Con1}\,(1)).
If $(i, j)= (m, n)\neq (s, t)\neq (1, 1)$,
then (\ref{f-3}) gives $[\vp^m\e^n, \vp^s\e^t ]=(s-t)\vp^s\e^t,$ which is the second case of (\ref{bracket-4-2}).
It remains to verify the last case of
(\ref{bracket-4-2}).  By (\ref{f-3}), we  have $[\vp^i\e^j, \vp^s\e^t] =0$ if $i, s \neq 1, m$
or $j, t \neq 1, n$.
We
 discuss the situations in two subcases.\vskip4pt

\noindent{\it Subcase 1.} Assume $i=1$ (thus $j\ne1$). Then  (\ref{f-3}) becomes
\begin{equation}\label{f-1-1}\begin{array}{lll}
[\vp\e^j, \vp^s\e^t]&\!\!\!\!=\!\!\!\!&\big(2n-m(j\!+\!1)\big)\d_{s,1}\d_{t,1}\vp^{2-m}\e^{j+1-n}+(s\!-\!t)\d_{1,m}\d_{j,n}\vp^s\e^t
-(1\!-\!j)\d_{s,m}\d_{t,n}\vp\e^j.
\end{array}
\end{equation}
Note that
 $\vp^{2-m}=0$ if $m>2$, 
in this case, (\ref{f-1-1}) becomes
$[\vp\e^j, \vp^s\e^t]=
-(1-j)\d_{s,m}\d_{t,n}\vp\e^j$,
and we have (\ref{bracket-4-2}).
Now assume $m=0$.  Then (\ref{f-1-1}) gives 
$[\vp\e^j, \vp^s\e^t]=2n\d_{s,1}\d_{t,1}\vp^{2}\e^{j+1-n}
-(1-j)\d_{s,0}\d_{t,n}\vp\e^j,$
and we see that
%
the last case of (\ref{bracket-4-2}) holds in this case.
Next assume $ m=1$. Then (\ref{f-1-1}) becomes 
%
%
$
[\vp\e^j, \vp^s\e^t]=(2n-(j+1))\d_{s,1}\d_{t,1}\vp\e^{j+1-n}+(s-t)\d_{j,n}\vp^s\e^t
-(1-j)\d_{s,1}\d_{t,n}\vp\e^j.
$ 
%
%
%
Hence the last case of (\ref{bracket-4-2}) holds. 
Finally assume $m=2$. By (\ref{f-1-1}), we have
\begin{equation}\label{f-1-3}\begin{array}{lll}
[\vp\e^j, \vp^s\e^t]&=&(2n-2(j+1))\d_{s,1}\d_{t,1}\e^{j+1-n}
-(1-j)\d_{s,2}\d_{t,n}\vp\e^j,\end{array}
\end{equation}
and the last case of (\ref{bracket-4-2}) holds again.
%
%
%
\vskip4pt

\noindent{\it Subcase 2.} Assume $i=m\neq 1 \neq s. $  We have
$[\vp^m\e^j, \vp^s\e^t]=(s-t)\d_{j,n}\vp^s\e^t-(m-j)\d_{s,m}\d_{t,n}\vp^m\e^j$ by (\ref{f-3}), i.e.,
$[\vp^m\e^j, \vp^s\e^t]=(s-t)\vp^s\e^t$ if $j=n$, or
$(j-m)\vp^m\e^j$ if $(s, t)=(m, n)$, or
$0$ otherwise.
This completes the proof of the theorem.
%
%
\QED

\begin{prop}\label{b-6-0} Let $f(x,y)= \sum
_{i=0}^{m}\sum
_{j=0}^n a_{i,j}x^iy^j \in \F[x,y]$ with $a_{m,n} \neq 0$, and
 $k,l\in\Z_+,$ $c\in\F$.
Denote ${\rm Supp\,}f=\{(i,j)\in\Z_+^2\,|\,a_{ij}\ne0\}$.
Then
\begin{eqnarray}\label{eq111}
{}\!\!\!\![x^ky^l, f(x, y)]= cf(x, y)\neq 0&\Longleftrightarrow&k=l=1,\ \ j-i=c,\ \forall\,(i,j)\in{\rm Supp\,}f,\\
\label{eq111+}
{}[x^ky^l, f(x, y)] =0\ \ \ \ \ \ \ \ \ \ \ \ \ &\Longleftrightarrow&kj-li=0, \ \forall\,(i,j)\in{\rm Supp\,}f.
\end{eqnarray}
In particular, if denote $r=A\otimes B- B\otimes A,$ where either $A=xy,$ $B= \sum
_{i=0}^n a_ix^iy^{c+i}$ for some $c \in \Z_+$ and $a_i\in\F$, or $ A=x^ky^l$, $B=f(x, y)=\sum
_{i=0}^{m}\sum
_{j=0}^n a_{i,j}x^iy^j$ such that $kj-li=0$ for all $(i, j)\in{\rm Supp\,}f $, then $r$ is a solution
of classical Yang-Baxter equation.
\end{prop}
\ni{\it Proof.~}~ We have
$[x^ky^l, f(x, y)]
\!=\!\sum_{i,j}a_{i,j}(kj-li)x^{k+i-1}y^{l+j-1}.$
By comparing the coefficients of the highest term (i.e., $x^my^n$) and $x^iy^j$ for all $i,j$, one immediately obtains \eqref{eq111} and \eqref{eq111+}.
%
%
%
\QED\vskip4pt

First consider $A=xy,\, B=\sum
_{i=0}^na_ix^iy^{i+m}$ with $a_n \neq 0$ and $m\in\Z_+\backslash\{0\}$. 
Then $[A, B]= mB\neq0$.
The triple $(\cal{A}, [\cdot, \cdot], \D_r)$ with $r= A\otimes B - B\otimes A$ is a coboundary triangular Lie bialgebra of Poisson type with bracket defined by
 (\ref{Poi4}) and cobracket 
 defined  by $$\D_r (g)
 = [g, A] \otimes B + A \otimes [g, B] - B \otimes [g, A] - [g, B] \otimes A,\ \forall\, g\in \AA.$$  The following  is one of the main results of the present paper.

 \begin{theo}\label{th-2-2} Let $(\cal{A}, [\cdot, \cdot], \D_r)$ be the coboundary triangular Lie bialgebra defined as above.
The dual Lie bialgebra of $(\cal{A}, [\cdot, \cdot], \D_r)$ is $(\cal{A}^{\circ}, [\cdot, \cdot], \D)$,  where  $\AA^{\circ}$ is described by Proposition $\ref{F[x]}(2),$
with cobracket $\D$ defined by \eqref{dual-1} and bracket uniquely determined by the skew-symmetry and  the following.
\begin{itemize}\item[\rm(1)]
In case $m\neq 0$, we have $($cf.~Convention $\ref{Con1}\,(1)$\,$)$, for $(s,t)\ne(1,1)$ and $p\ne q$, 
\begin{eqnarray}
\label{eq2.2-1}
&\!\!\!\!\!\!\!\!\!\!\!\!\!\!\!\!\!\!\!\!\!\!\!\!\!&
[\vp^p\e^p, \vp^s\e^t]=\d_{p,1} \mbox{$\sum\limits_{i=0}^n$} a_i(si\!+\!sm\!+\!m\!-\!ti)\vp^{s-i+1} \e^{t-i-m+1}
,\\
\label{eq2.2-2}
&\!\!\!\!\!\!\!\!\!\!\!\!\!\!\!\!\!\!\!\!\!\!\!\!\!&
[\vp^p\e^q, \vp^s\e^t]=(p\!-\!q)\mbox{$\sum\limits_{i=0}^n$} a_i\d_{s,i}\d_{t, i+m}\vp^p\e^q -(s\!-\!t)\mbox{$\sum\limits_{i=0}^n$} a_i\d_{p,i}\d_{q, i+m}\vp^s\e^t
.
\end{eqnarray}
\item[\rm(2)] In case $m=0$, we have, for $(s,t)\ne(1,1)$,
\begin{eqnarray}
\label{+eq2.2-1}
\!\!\!\!\!\!\!\!\!\!\!\!\!\!\!\!\!\!\!\!\!\!\!\!&\!\!\!\!\!\!\!\!\!\!\!\!\!\!\!\!\!\!\!\!\!\!\!\!\!&
 [\vp\e, \vp^s\e^t]  =\mbox{$\sum\limits_{i=2}^n$}a_i (s-t)i\vp^{s-i+1} \e^{t-i+1}
,\\
\label{+eq2.2-2}
\!\!\!\!\!\!\!\!\!\!\!\!\!\!\!\!\!\!\!\!\!\!\!\!&\!\!\!\!\!\!\!\!\!\!\!\!\!\!\!\!\!\!\!\!\!\!\!\!\!&
[\vp^p\e^p, \vp^s\e^t] =(t-s)a_p\vp^s\e^t \mbox{ \ if \ }\ p\in \{0, 2, 3, \cdots, n\},\\
\label{+eq2.2-3}
\!\!\!\!\!\!\!\!\!\!\!\!\!\!\!\!\!\!\!\!\!\!\!\!&\!\!\!\!\!\!\!\!\!\!\!\!\!\!\!\!\!\!\!\!\!\!\!\!\!&
 [\vp^p\e^q, \vp^s\e^t] = 0 \mbox{ \ if \ }\ p\neq q,\, s\neq t.
 \end{eqnarray}
\end{itemize}
\end{theo}

\ni{\it Proof.~}~ Denote $C\!=\![x^ky^l, xy]\!=\!(k\!-\!l)x^ky^l$ and
$D\! =\![x^ky^l, B]
\!=\!\mbox{$ \sum
_{i=0}^n$}a_i(ki\!+\!km\!-\!li)x^{i+k-1}y^{i+l+m-1}.$
We have
\begin{eqnarray}\label{theorem-2-3}
&\!\!\!\!\!\!\!\!\!\!\!\!\!\!\!\!\!\!\!\!\!&
\langle[\vp^p\e^q, \vp^s\e^t], x^ky^l\rangle
=\langle\vp^p\e^q\otimes \vp^s\e^t, C \otimes B + A\otimes D- B\otimes C - D \otimes A\rangle
=P^{p,q}_{s,t} - P^{s,t}_{p,q}
%
%
,\end{eqnarray}
where (regarding $k,l$ as fixed)
$$\begin{array}{lll} P^{p,q}_{s,t}&\!\!\!=\!\!\!&(k-l)\d_{p,k}\d_{q,l}\sum\limits_{i=0}^n a_i\d_{s,i}\d_{t, i+m}
 + \d_{p,1}\d_{q,1}\sum\limits_{i=0}^n a_i(ki+km-li)\d_{s, i+k-1}\d_{t, i+m+l-1}\\[12pt]
&\!\!\!=\!\!\!& (p-q)\d_{p,k}\d_{q,l}\sum\limits_{i=0}^n a_i\d_{s,i}\d_{t, i+m} +
 \d_{p,1}\d_{q,1}\sum\limits_{i=0}^na_i(si+sm+m-ti)\d_{s-i+1, k}\d_{t-i-m+1,l}
 \\[12pt]
&\!\!\!=\!\!\!&\langle H^{p,q}_{s,t}, x^ky^l\rangle,\mbox{ \ and where,}\end{array}$$
\begin{equation}\label{theorem-2-4}H^{p,q}_{s,t} = (p-q)\mbox{$\sum\limits_{i=0}^n$} a_i\d_{s,i}\d_{t, i+m}\vp^p\e^q + \d_{p,1}\d_{q,1}
\mbox{$\sum\limits_{i=0}^n$}a_i (si+sm+m-ti)\vp^{s-i+1} \e^{t-i-m+1}.\end{equation}
Thus
\begin{equation}\label{theorem-2-6}
[\vp^p\e^q, \vp^s\e^t] = H^{p,q}_{s,t} - H^{s,t}_{p,q}.
\end{equation}Assume $(s,t)\ne(1,1)$.
First suppose $m\ne0$.
%
%
Then \eqref{theorem-2-3}--\eqref{theorem-2-6} give $[\vp^p\e^p, \vp^s\e^t] =0$ for $p=q\ne1$, and
%
$[\vp\e, \vp^s\e^t]=\sum_{i=0}^n a_i(si+sm+m-ti)\vp^{s-i+1} \e^{t-i-m+1}
$ for $p=q=1$.
%
%
We have \eqref{eq2.2-1}.
%
%
If $p\neq q$, we have
$ 
[\vp^p\e^q, \vp^s\e^t]=
(p-q)\sum
_{i=0}^n a_i\d_{s,i}\d_{t, i+m}\vp^p\e^q -(s-t)\sum
_{i=0}^n a_i\d_{p,i}\d_{q, i+m}\vp^s\e^t
,$ 
by \eqref{theorem-2-3}--\eqref{theorem-2-6},
and we have \eqref{eq2.2-2}.
Now suppose $m=0$. 
%
Then 
\eqref{theorem-2-4} becomes
\begin{equation}\label{theorem-2-7}H^{p,q}_{s,t} = (p-q)\mbox{$\sum\limits_{i=0}^n$} a_i\d_{s,i}\d_{t, i}\vp^p\e^q +
 \d_{p,1}\d_{q,1}\mbox{$\sum\limits_{i=0}^n$}a_i (s-t)i\vp^{s-i+1} \e^{t-i+1}.\end{equation}
%
We have
$[\vp\e, \vp^s\e^t] =\sum
_{i=2}^na_i (s-t)i\vp^{s-i+1} \e^{t-i+1}$ for $(p,q)=(1, 1)$, 
i.e., we have \eqref{+eq2.2-1}.
%
If  $p=q\in\{0,2,3,...,n\}$, 
it is easy to see from \eqref{theorem-2-7} and \eqref{theorem-2-6} that we have \eqref{+eq2.2-2}.
%
%
%
%
Finally assume $p\neq q$. One can easily obtain \eqref{+eq2.2-3}.
\QED

\begin{rema}\rm\label{Raa22} Theorems \ref{theorem-2} and \ref{th-2-2} and the following theorem provide us some examples of
the nontriviality of coboundary triangular  Lie bialgebras, even in the case of very trivial solution $r=A\otimes B-B\otimes A$ of CYBE
with $[A, B]=0$.
\end{rema}

\begin{theo}\label{theorem-3-1} Let $(k,l)\in\Z_+^2$ be fixed and $A=x^ky^l, B=f(x, y)=\sum_{(i, j) \in S}a_{i, j}x^iy^j \in \AA$ with $a_{i,j}\in\F$,  where
$S={\rm Supp\,}f$ is some subset of $\Z_+^2$ such that $kj-li=0$ for $(i,j)\in S$. Denote $r=A\otimes B- B\otimes A$.
Then $(\cal{A}, [\cdot, \cdot], \D_r)$ is a coboundary triangular Lie bialgebra
of Poisson type. The dual Lie bialgebra of $(\cal{A}, [\cdot, \cdot], \D_r)$ is $(\cal{A^{\circ}}, [\cdot, \cdot], \D)$ with cobracket $\D$ defined as in Theorem $\ref{theorem-2}$ and bracket uniquely determined by the skew-symmetry and the following.\vskip4pt

$(1)$ If  $(s, t) \neq (k, l)$, then
{
\begin{eqnarray}\label{caseqq1}
[\vp^k\e^l, \vp^s\e^t]=
\left\{\begin{array}{ll}\!\!
(l\!-\!k)a_{s, t}\vp \e\!+\!\sum\limits_{(i,j)\in S}a_{i, j}(j\!-\!i)\vp^{s-i+1}\e^{t-j+1}
\!-\!(l\!-\!k)a_{k,l}\vp^{s-k+1}\e^{t-l+1}
,\\[12pt]\!\!\!
\sum\limits_{(i,j)\in S}\!\!a_{i, j}(j(s\!+\!1)\!-\!i(t\!+\!1))\vp^{s-i+1}\e^{t-j+1}\!-\!(l(s\!+\!1)\!-\!k(t\!+\!1))a_{k,l}\vp^{s-k+1}\e^{t-l+1}
,\\[12pt]
(l\!-\!k)a_{s, t}\vp \e \!+\!\sum\limits_{(i,j)\in S}a_{i, j}(j\!-\!i)\vp^{s-i+1}\e^{t-j+1}
,\\[12pt]
\sum\limits_{(i,j)\in S}a_{i, j}(j(s+1)-i(t+1))\vp^{s-i+1}\e^{t-j+1}
,\end{array}\right.
\end{eqnarray}%
}%
according to the following four cases
$$\mbox{{\rm(i)} $(s, t), (k, l) \!\in\! S$, {\rm(ii)} $(s, t) \!\notin \!S,(k,l)\! \in\! S$,
{\rm(iii)} $(s, t)\!\in\! S, (k,l)\! \notin\! S$ or {\rm (iv)} $(s, t),(k, l)\notin S$.}$$

$(2)$ If $(p, q)\neq (k, l),\, (s, t) \neq (k, l)$, then
\begin{eqnarray}\label{caseqq2}
[\vp^p\e^q, \vp^s\e^t ]=\left\{\begin{array}{ll}
(l\!-\!k)a_{s, t}\vp^{p-k+1}\e^{q-l+1}-(l\!-\!k)a_{p, q}\vp^{s-k+1}\e^{t-l+1}\!\!\!&\mbox{if}\ (p, q), (s, t) \in S,\\[4pt]
(l(p\!+\!1)\!-\!k(q\!+\!1))a_{s, t}\vp^{p-k+1}\e^{q-l+1}& \mbox{if}\ (p, q)\!\notin\! S, (s, t)\!\in\! S,\\[4pt]
0&\mbox{if}\ (p, q)\notin S, (s, t) \notin S.\end{array}\right.\end{eqnarray}
\end{theo}

\ni{\it Proof.~}~ Since $[A, B] =0,$ the triple $(\cal{A}, [\cdot, \cdot], \D_r)$ is obviously a coboundary triangular Lie bialgebra.  We only need to determine the bracket relations.
Note that $\langle [\vp^p\e^q, \vp^s\e^t], x^{m^{\prime}}y^{n^{\prime}} \rangle$ equals
\begin{eqnarray}\label{theorem-3-4}
&\!\!\!\!\!\!\!\!\!\!\!\!&
\phantom{X=}\Big\langle \vp^p\e^q\otimes\vp^s\e^t, [x^{m^{\prime}}y^{n^{\prime}}, A] \otimes B + A\otimes [x^{m^{\prime}}y^{n^{\prime}}, B]
-B\otimes[x^{m^{\prime}}y^{n^{\prime}}, A] -[x^{m^{\prime}}y^{n^{\prime}}, B]\otimes A \Big\rangle
\nonumber\\&\!\!\!\!\!\!\!\!\!\!\!\!&
\phantom{X}
=
(m^{\prime}l-n^{\prime}k)\d_{p,{m^{\prime}+k-1}}\d_{q,{n^{\prime}+l-1}}\mbox{$\sum\limits_{(i, j)\in S}$}a_{i, j}\d_{s,i}\d_{t,j}
+\d_{p,k}\d_{q, l}\mbox{$\sum\limits_{(i,j)\in S}$}a_{i, j}(m^{\prime}j - n^{\prime} i)\d_{s,{i+m^{\prime}-1}}\d_{t, {j+n^{\prime}-1}}
\nonumber\\&\!\!\!\!\!\!\!\!\!\!\!\!&
\phantom{X=}
-(m^{\prime}l\!-\!n^{\prime}k)\d_{s, {m^{\prime}+k-1}}\d_{t,{n^{\prime}+l-1}}\mbox{$\sum\limits_{(i, j)\in S}$}a_{i, j}\d_{p,i}\d_{q,j}
\!-\!\d_{s,k}\d_{t, l}\mbox{$\sum\limits_{(i, j)\in S}$}a_{i, j}(m^{\prime}j\! -\! n^{\prime} i)\d_{p,{i+m^{\prime}-1}}\d_{q, {j+n^{\prime}-1}}
%
\nonumber\\&\!\!\!\!\!\!\!\!\!\!\!\!&
\phantom{X}
=\langle H^{p, q}_{s, t}- H^{s, t}_{p, q}, x^{m^{\prime}}y^{n^{\prime}}\rangle ,
\end{eqnarray}
where (noting that $a_{s,t}=0$ if $(s,t)\notin S$)
\begin{equation*}
\begin{array}{lll}H^{p, q}_{s, t}\!=\!\big(l(p\!+\!1)\!-\!k(q\!+\!1)\big)a_{s,t}\vp^{p-k+1}\e^{q-l+1}
\!+\!\d_{p,k}\d_{q, l}\sum\limits_{(i,j)\in S}a_{i, j}\big(j(s\!+\!1)\!-\!i(t\!+\!1)\big)\vp^{s-i+1}\e^{t-j+1}.\end{array}\end{equation*}
%
%
If  $(p, q)= (k, l)$ and $(s, t) \neq (k, l)$, then 
\eqref{theorem-3-4} gives (noting that $a_{k,l}=0$ if $(k,l)\notin S$)
$$\begin{array}{lll}[\vp^k\e^l, \vp^s\e^t]& =&(l-k)a_{s,t}\vp\e 
 +\sum\limits_{(i,j)\in S}a_{i, j}\big(j(s+1)-i(t+1)\big)\vp^{s-i+1}\e^{t-j+1}\\[12pt]
&&- \big(l(s+1)-k(t+1)\big)a_{k,l}\vp^{s-k+1}\e^{t-l+1}
.\end{array}$$
In particular, if $(s, t), (k, l) \in S$, then (using the fact that $kj-li=0$ for  $(i,j)\in S$)
$$[\vp^k\e^l, \vp^s\e^t]=(l-k)a_{s, t}\vp \e +\mbox{$\sum\limits_{(i,j)\in S}$}a_{i, j}(j-i)\vp^{s-i+1}\e^{t-j+1}
-(l-k)a_{k,l}\vp^{s-k+1}\e^{t-l+1},$$
which gives the first case of \eqref{caseqq1}.
If $(s, t) \notin S, (k,l) \in S$, then
$$[\vp^k\e^l, \vp^s\e^t]=\mbox{$\sum\limits_{(i,j)\in S}$}a_{i, j}\big(j(s\!+\!1)-i(t\!+\!1)\big)\vp^{s-i+1}\e^{t-j+1}-\big(l(s\!+\!1)-k(t\!+\!1)\big)a_{k,l}\vp^{s-k+1}\e^{t-l+1},$$
which gives the second case of \eqref{caseqq1}.
If $(s, t)\in S, (k,l) \notin S$, then
$[\vp^k\e^l, \vp^s\e^t]=(l-k)a_{s, t}\vp \e +\mbox{$\sum
_{(i,j)\in S}$}a_{i, j}(j-i)\vp^{s-i+1}\e^{t-j+1},$
which gives the third case of \eqref{caseqq1}.
If $(s, t) \notin S, (k, l)\notin S$, then
$[\vp^k\e^l, \vp^s\e^t]=\mbox{$\sum
_{(i,j)\in S}$}a_{i, j}(j(s+1)-i(t+1))\vp^{s-i+1}\e^{t-j+1},$
which completes the proof of \eqref{caseqq1}.

Now assume  $(p, q)\neq (k, l), (s, t) \neq (k, l)$. Then  
 \eqref{theorem-3-4} gives
$$\begin{array}{lll}[\vp^p\e^q, \vp^s\e^t ]=
 \big(l(p+1)-k(q+1)\big)a_{s,t}\vp^{p-k+1}\e^{q-l+1}
-\big(l(s+1)-k(t+1)\big)a_{p,q}\vp^{s-k+1}\e^{t-l+1}
.\end{array}$$
If $(p, q), (s, t) \!\in\! S$, then
$[\vp^p\e^q, \vp^s\e^t ]\!=\!(l\!-\!k)a_{s, t}\vp^{p-k+1}\e^{q-l+1}\!-\!(l\!-\!k)a_{p, q}\vp^{s-k+1}\e^{t-l+1},$
giving the first case of \eqref{caseqq2}. If
$(p, q)\!\notin\! S, (s, t)\!\in\! S$, then
$[\vp^p\e^q, \vp^s\e^t ]\!=\!\big(l(p\!+\!1)\!-\!k(q\!+\!1)\big)a_{s, t}\vp^{p-k+1}\e^{q-l+1},$
giving the second case of \eqref{caseqq2}. In case
%
$(p, q)\notin S, (s, t) \notin S$, we have
$[\vp^p\e^q, \vp^s\e^t] =0$, which completes the proof of the theorem.
\QED\vskip4pt

In the final part of the paper, we will present an example of a dual Lie bialgebra which has different feature from the previous dual Lie bialgebras (Theorems \ref{theorem-2}, \ref{th-2-2} and \ref{theorem-3-1}).

Denote $A=x(1+y)(2+y),$ $B=x^2(1+y)^3(2+y) \in \AA$. Then \eqref{Poi4} shows  $[A, B] =B$. Take
 $r= A\otimes B- B\otimes A$. We obtain  a coboundary triangular Lie bialgebra
$(\AA, [\cdot, \cdot], \D)$ with bracket defined by (\ref{Poi4}), and  cobracket defined by
$$\D(f) = f\cdot r=[f,A]\otimes B+A\otimes [f,B]-[f,B]\otimes A-B\otimes [f,A].$$
\begin{theo}\adddot \label{last-theo}Let $(\AA, [\cdot, \cdot], \D)$ be the Lie bialgebra defined above. The dual Lie bialgebra
of $(\AA, [\cdot, \cdot], \D)$ is $(\AA^{\circ}, [\cdot, \cdot], \D)$, where  $\AA^{\circ}$ and $\D$ are
 defined in Theorem $\ref{theorem-2}$, and the bracket is uniquely determined by the skew-symmetry and the following.
 First, denote
 $$\begin{array}{ll} c_0=2,\  c_1=3, \ c_2=1, \mbox{ \ and \ }  c_j=0\mbox{ for }j\neq 0,1,2,\\[2pt]
k_0=2,\ k_1=7,\ k_2=9,\ k_3=5,\ k_4=1, \mbox{ \ and \ } k_t=0 \mbox{ for }\ 4<t \in \Z_+.\end{array}
$$
Then
\begin{eqnarray}
\label{b-5-1}&\!\!\!\!\!\!\!\!\!\!\!\!\!\!\!\!\!\!\!\!\!\!\!\!\!\!\!\!\!\!&
 \begin{array}{lll}  [\vp\e^j, \vp^s\e^t]&\!\!\!=\!\!\!&c_j\Big((4s\!-\!2t\!+\!2)\vp^{s-1}\e^{t-3}\!+\!(15s\!-\!10t\!+\!5)\vp^{s-1}\e^{t-2}\\[4pt] &&
 +18(s\!-\!t)\vp^{s-1}\e^{t-1}
 \!+\!(7s\!-\!14t\!-\!7)\vp^{s-1}\e^t
 \!+\!4(t\!+\!1)\vp^{s-1}\e^{t+1}\Big)\mbox{ if }s\ne1,2, 
\end{array}
\end{eqnarray}

\begin{eqnarray}
\label{b-5-2}&\!\!\!\!\!\!\!\!\!\!\!\!\!\!\!\!\!\!\!\!\!\!\!\!\!\!\!\!\!\!&
\begin{array}{lll}
[\vp\e^j, \vp\e^t] &\!\!\!=\!\!\!&c_j\Big((6\!-\!2t)\e^{t-3}
\!+\!(20\!-\!10t)\e^{t-2}\! +\!18(1\!-\!t)\e^{t-1}
\!-\!14t\e^t \!+\!4(t\!+\!1)\e^{t+1}\Big)\\
&&-c_t\Big((6\!-\!2j)\e^{j-3}\!+\!(20\!-\!10j)\e^{j-2} 
\!+\!18(1\!-\!j)\e^{j-1}\!-\!14j\e^j\! +\!4(j\!+\!1)\e^{j+1}\Big),
\end{array}\\[-5pt]
\label{b-5-3}
&\!\!\!\!\!\!\!\!\!\!\!\!\!\!\!\!\!\!\!\!\!\!\!\!\!\!\!\!\!\!&
\begin{array}{lll}
[\vp\e^j, \vp^2\e^t] &\!\!\!=\!\!\!&k_t\Big((3-j)\vp\e^{j-1}+3(1-j)\vp\e^j-2(j+1)\vp\e^{j+1}\Big)\\
 &&+ c_j\Big((10\!-\!2t)\vp\e^{t-3}
%
%
%
\\
&&\phantom{XX} \!+\!(35\!-\!10t)\vp\e^{t-2} \!+\!18(2\!-\!t)\vp\e^{t-1}\!+\!(7\!-\!14t)\vp\e^t \!+\!4(t\!+\!1)\vp\e^{t+1}\Big),
\end{array}
\\[-5pt]
\label{b-5-4}&\!\!\!\!\!\!\!\!\!\!\!\!\!\!\!\!\!\!\!\!\!\!\!\!\!\!\!\!\!\!&
\begin{array}{lll}
[\vp^2\e^j, \vp^s\e^t] = -k_j\Big((2s-t+1)\vp^s\e^{t-1}+3(s-t)\vp^s\e^t-2(t+1)\vp^s\e^{t+1}\Big)\mbox{ if }s\neq 1,2,
\end{array}
\\[-5pt]
\label{b-5-5}
&\!\!\!\!\!\!\!\!\!\!\!\!\!\!\!\!\!\!\!\!\!\!\!\!\!\!\!\!\!\!&
\begin{array}{lll}
[\vp^2\e^j, \vp^2\e^t]&\!\!\!=\!\!\!&k_t\Big((5-j)\vp^2\e^{j-1}+3(2-j)\vp^2\e^j-2(j+1)\vp^2\e^{j+1}\Big) \\
 &&-k_j\Big((5-t)\vp^2\e^{t-1}+3(2-t)\vp^2\e^t-2(t+1)\vp^2\e^{t+1}\Big),\end{array}
\\[-5pt]
\label{b-5-6}
&\!\!\!\!\!\!\!\!\!\!\!\!\!\!\!\!\!\!\!\!\!\!\!\!\!\!\!\!\!\!&
\begin{array}{lll}
[\vp^i\e^j, \vp^s\e^t ]= 0 \mbox{ \ if \ }\ i\neq 1,2, s\neq 1,2.
\end{array}
\end{eqnarray}
\end{theo}
\ni{\it Proof.~}~ Denote $C= [x^my^n, A], D=[x^my^n, B]$. By (\ref{Poi4}), we \vspace*{-5pt}have
\begin{eqnarray*}
&\!\!\!\!\!\!\!\!\!\!\!\!\!\!\!\!\!\!\!\!\!\!\!\!\!\!\!\!\!\!&
C= [x^my^n, x(y+1)(y+2)]= (2m-n)x^my^{n+1} + 3(m-n)x^my^n -2nx^my^{n-1},\\
&\!\!\!\!\!\!\!\!\!\!\!\!\!\!\!\!\!\!\!\!\!\!\!\!\!\!\!\!\!\!&
D=[x^my^n, x^2(y+1)^3(y+2)]
\\[-2pt]
&\!\!\!\!\!\!\!\!\!\!\!\!\!\!\!\!\!\!\!\!\!\!\!\!\!\!\!\!\!\!&
\phantom{D}
=(4m-2n)x^{m+1}y^{n+3} + (15m-10n)x^{m+1}y^{n+2}
+18(m-n)x^{m+1}y^{n+1}
\\[-2pt]
&\!\!\!\!\!\!\!\!\!\!\!\!\!\!\!\!\!\!\!\!\!\!\!\!\!\!\!\!\!\!&
\phantom{D=}
 +(7m-14n)x^{m+1}y^n -4nx^{m+1}y^{n-1}.\end{eqnarray*}
For $i, j, s, t \in \Z_+$, note that
$\langle [\vp^i\e^j, \vp^s\e^t], x^my^n\rangle$ is equal to
\begin{eqnarray}
&\!\!\!\!\!\!\!\!\!\!\!\!&
\phantom{XX=}
\langle \vp^i\e^j \otimes \vp^s\e^t, [x^my^n, A]\otimes B + A\otimes [x^my^n, B]
-[x^my^n, B] \otimes A -B \otimes [x^my^n, A] \rangle
\nonumber\\
&\!\!\!\!\!\!\!\!\!\!\!\!&
\phantom{XX}
=\langle \vp^i\e^j \otimes \vp^s\e^t, C\otimes B + A\otimes D- D\otimes A-B\otimes C \rangle
=
P^{i,j}_{s,t} - P^{s,t}_{i,j},
\end{eqnarray}
where
$ P^{i,j}_{s,t}= \langle \vp^i\e^j, C \rangle \langle \vp^s\e^t, B \rangle + \langle \vp^i\e^j, A \rangle \langle \vp^s\e^t, D \rangle$, which is equal \vspace*{-5pt}to
$$\begin{array}{lll}
 &\!\!\!\phantom{=}\!\!\!& \Big((2m\!-\!n)\d_{i,m}\d_{j, n+1}\!+\!3(m\!-\!n)\d_{i,m}\d_{j,n}\!-\!2n\d_{i,m}\d_{j,n-1}\Big)\\[6pt]
&\!\!\!\!\!\!&\phantom{\Big(}\times\Big(2\d_{s,2}\d_{t,0}+7\d_{s,2}\d_{t,1}+9\d_{s,2}\d_{t,2}+5\d_{s,2}\d_{t,3}+\d_{s,2}\d_{t,4}\Big)
\\[4pt]
&\!\!\!\!\!\!&+\Big(2\d_{i,1}\d_{j,0} + 3\d_{i,1}\d_{j,1}+\d_{i,1}\d_{j,2}\Big)\\[4pt]
&\!\!\!\!\!\!&\phantom{\Big(}\times\Big((4m\!-\!2n)\d_{s,m+1}\d_{t,n+3}\!+\! (15m\!-\!10n)\d_{s,m+1}\d_{t,n+2} \!+\!18(m\!-\!n)\d_{s,m+1}\d_{t,n+1}\\[4pt]
 &\!\!\!\!\!\!&\phantom{\Big(\times\Big(}+ (7m-14n)\d_{s, m+1}\d_{t,n}-4n\d_{s,m+1}\d_{t,n-1}\Big)
\\[6pt]
&\!\!\!=\!\!\!&\Big\langle \Big(2\d_{s,2}\d_{t,0}+7\d_{s,2}\d_{t,1}+9\d_{s,2}\d_{t,2}+5\d_{s,2}\d_{t,3}+\d_{s,2}\d_{t,4}\Big)\\[4pt]
&\!\!\!\!\!\!&\phantom{\Big\langle+}\times\Big((2i-j+1)\vp^i\e^{j-1}+3(i-j)\vp^i\e^j-2(j+1)\vp^i\e^{j+1})  \\[4pt]
&\!\!\!\!\!\!&\phantom{\Big\langle}+\Big(2\d_{i,1}\d_{j,0} + 3\d_{i,1}\d_{j,1}+\d_{i,1}\d_{j,2}\Big)\\[4pt]
&\!\!\!\!\!\!&\phantom{\Big\langle+}\times\Big((4s-2t+2)\vp^{s-1}\e^{t-3}
+(15s-10t+5)\vp^{s-1}\e^{t-2} +18(s-t)\vp^{s-1}\e^{t-1}\\[4pt]
 &\!\!\!\!\!\!&\phantom{\Big\langle+\times}+(7s-14t-7)\vp^{s-1}\e^t +4(t+1)\vp^{s-1}\e^{t+1}\Big),x^my^n\Big\rangle.\end{array}$$
We obtain
\begin{equation}\label{H-12}[\vp^i\e^j, \vp^s\e^t ]= H^{i,j}_{s,t} - H^{s, t}_{i, j},\end{equation}
where $H^{i,j}_{s,t}$ is equal to \begin{eqnarray*}
&\!\!\!\!\!\!&
\Big(2\d_{s,2}\d_{t,0}\!+\!7\d_{s,2}\d_{t,1}\!+\!9\d_{s,2}\d_{t,2}\!+\!5\d_{s,2}\d_{t,3}\!+\!\d_{s,2}\d_{t,4}\Big)
\Big((2i\!-\!j\!+\!1)\vp^i\e^{j-1}\!+\!3(i\!-\!j)\vp^i\e^j\!-\!2(j\!+\!1)\vp^i\e^{j+1}\Big) \nonumber\\[-4pt]
&\!\!\!\!\!\!&+\Big(2\d_{i,1}\d_{j,0} + 3\d_{i,1}\d_{j,1}+\d_{i,1}\d_{j,2}\Big)\nonumber\\[-4pt]
&\!\!\!\!\!\!&\phantom{\Big(}\times\Big((4s\!-\!2t\!+\!2)\vp^{s-1}\e^{t-3}\!+\!(15s\!-\!10t\!+\!5)\vp^{s-1}\e^{t-2} \!+\!18(s\!-\!t)\vp^{s-1}\e^{t-1}
\nonumber\\[-4pt]
&\!\!\!\!\!\!&\phantom{\Big(\times\Big(}+(7s-14t-7)\vp^{s-1}\e^t +4(t+1)\vp^{s-1}\e^{t+1}\Big).\nonumber\end{eqnarray*}
%
%
If  $i \!\neq\! 1,2$, $s\!\neq\! 1,2$,
then 
\eqref{H-12} gives that $[\vp^i\e^j, \vp^s\e^t ]=0$, which is \eqref{b-5-6}.
Assume $i=1,\, s\neq 1,2$. Then
$[\vp\e^j, \vp^s\e^t] = H^{1, j}_{s,t}$, which gives
\eqref{b-5-1}.
If $i=1, s=1,$ then 
$[\vp\e^j, \vp\e^t] $ is equal to
$$\begin{array}{lll}\Big(2\d_{j,0} + 3\d_{j,1}+\d_{j,2}\Big)\Big((6-2t)\e^{t-3}
+(20-10t)\e^{t-2}
 +18(1-t)\e^{t-1}-14t\e^t +4(t+1)\e^{t+1}\Big)\\[4pt]
-\Big(2\d_{t,0} \!+\! 3\d_{t,1}\!+\!\d_{t,2}\Big)\Big((6\!-\!2j)\e^{j-3}\!+\!(20\!-\!10j)\e^{j-2}
  \!+\!18(1\!-\!j)\e^{j-1}\!-\!14j\e^j \!+\!4(j\!+\!1)\e^{j+1}\Big),\end{array}$$
which implies (\ref{b-5-2}).
Assume $i=1, s=2$. Then 
(\ref{H-12}) implies that $[\vp\e^j, \vp^2\e^t]$ is equal to
$$\begin{array}{lll}\Big(2\d_{t,0}+7\d_{t,1}+9\d_{t,2}
+5\d_{t,3}+\d_{t,4}\Big)\Big((3-j)\vp\e^{j-1}+3(1-j)\vp\e^j-2(j+1)\vp\e^{j+1}\Big) \\[4pt]
+\Big(2\d_{j,0}\! +\! 3\d_{j,1}\!+\!\d_{j,2}\Big)\Big((10\!-\!2t)\vp\e^{t-3}\!+\!(35\!-\!10t)\vp\e^{t-2}\! +\!18(2\!-\!t)\vp\e^{t-1}\\[4pt]
\phantom{+\Big(2\d_{j,0}\! +\! 3\d_{j,1}\!+\!\d_{j,2}\Big)\Big(}+(7-14t)\vp\e^t +4(t+1)\vp\e^{t+1}\Big),\end{array}$$
which gives (\ref{b-5-3}).
If $i=2, s\neq 1,2$, then 
(\ref{H-12}) shows that
$[\vp^2\e^j, \vp^s\e^t]$ is equal to
$$\begin{array}{lll}\!\!\!\!-(2\d_{j,0}+7\d_{j,1}+9\d_{j,2}
+5\d_{j,3}+\d_{j,4})\Big((2s-t+1)\vp^s\e^{t-1}+3(s-t)\vp^s\e^t-2(t+1)\vp^s\e^{t+1}\Big),\end{array}$$
which is (\ref{b-5-4}).
Finally if $i=2, s=2$, we see from 
(\ref{H-12}) that $[\vp^2\e^j, \vp^2\e^t]$ is equal to
$$\begin{array}{lll}\Big(2\d_{t,0}+7\d_{t,1}+9\d_{t,2}+5\d_{t,3}+\d_{t,4}\Big)
\Big((5-j)\vp^2\e^{j-1}+3(2-j)\vp^2\e^j-2(j+1)\vp^2\e^{j+1}\Big) \\[4pt]
 -\Big(2\d_{j,0}+7\d_{j,1}+9\d_{j,2}+5\d_{j,3}+\d_{j,4}\Big)\Big((5-t)\vp^2\e^{t-1}+3(2-t)\vp^2\e^t-2(t+1)\vp^2\e^{t+1}\Big),\end{array}$$
and we obtain (\ref{b-5-5}). The proof of the theorem is completed.
\QED
\section{Conclusion remark}
Theorems \ref{theorem-2}, \ref{th-2-2}, \ref{theorem-3-1} and \ref{last-theo} provide 
five classes of infinite dimensional  Lie Bialgebras $(\AA^{\circ}, [\cdot, \cdot], \D)$ of Poisson type.
As by-products, we  obtain five new classes of infinite dimensional Lie algebras $({\cal L}, [\cdot, \cdot])$ with the underlining space ${\cal L}=\F[\vp,\e]$
(the polynomial algebra on two variables $\vp,\,\e$) 
and brackets defined respectively by (\ref{bracket-4-2}), (\ref{eq2.2-1})--(\ref{eq2.2-2}), (\ref{+eq2.2-1})--(\ref{+eq2.2-3}), (\ref{caseqq1})--(\ref{caseqq2}) and (\ref{b-5-1})--(\ref{b-5-6}). We close the paper by proposing the following questions: In which conditions will
 Lie bialgebras $(\AA^{\circ}, [\cdot, \cdot], \D)$ defined in Theorems \ref{theorem-2}, \ref{th-2-2}, \ref{theorem-3-1} and \ref{last-theo} be coboundary triangular? What
kinds of  structure and representation theories will these Lie algebras $({\cal L}, [\cdot, \cdot])$ have?

\Acknowledgements{ Supported by NSF grant  11071147, 11431010, 11371278
 of China,
NSF grant ZR2010AM003,  ZR2013AL013 of Shandong Province,
SMSTC grant 12XD1405000,  Fundamental Research Funds for the Central Universities.}


\end{document}